\documentclass[12pt,reqno]{article}

\usepackage[usenames]{color}
\usepackage{amssymb}
\usepackage{graphicx}
\usepackage{amscd}

\usepackage[colorlinks=true,
linkcolor=webgreen,
filecolor=webbrown,
citecolor=webgreen]{hyperref}

\definecolor{webgreen}{rgb}{0,.5,0}
\definecolor{webbrown}{rgb}{.6,0,0}

\usepackage{color}
\usepackage{fullpage}
\usepackage{float}
\usepackage{epsfig}

\usepackage{graphics,amsmath,amssymb}
\usepackage{amsthm}
\usepackage{amsfonts}
\usepackage{latexsym}
\usepackage{epsf}

\newcommand{\seqnum}[1]{\href{http://oeis.org/#1}{\underline{#1}}}

\begin{document}

\begin{center}
\epsfxsize=4in
%\leavevmode\epsffile{logo129.eps}
\end{center}

\begin{center}
\vskip 0.1cm{\LARGE\bf A Generalized Ap\'{e}ry Series}
\vskip 0.5cm
\large
M. L. Glasser \\
Department of Physics,\\
Clarkson University, Postdam, NY 13699-5820\\
\href{mailto:laryg@clarkson.edu}{\it laryg@clarkson.edu}\\

\end{center}

\vskip .2 in

\begin{abstract}
The inverse central binomial series
\begin{equation}S_k(z)=\sum_{n=1}^{\infty}\frac{n^k z^n}{\binom{2n}{n}}\nonumber\end{equation}
popularized by Ap\'ery and Lehmer is evaluated for positive 
integers $k$ along with the asymptotic behavior for large $k$. It is found 
that value $z=2$, as commented on by D.  H.  Lehmer provides a unique 
relation to $\pi$.

\end{abstract}

\bigskip

\section{Introduction}

Since the appearance of $S_{-3}(1)$ in Ap\'ery's famous proof \cite{aper} in 
1979 that $\zeta(3)$ is irrational, an extensive literature has been 
devoted to the series
\begin{equation}
S_k(z)=\sum_{n=1}^{\infty}\frac{n^k z^n}{\binom{2n}{n}}
\label{1}
\end{equation}
For example, in 1985   Lehmer \cite{lehm} presented a number of special cases 
which could be obtained from the Taylor series for 
$f(x)=x^{-1/2}(1-x)^{-1/2}\sin^{-1}x$ using only elementary calculus. In 
passing, he noted that when $k$ is a positive integer, $S_k(2)$ had the 
form $a_k-b_k\pi$ and that  the rational number $a_k/b_k$ ``is a close 
approximation to $\pi$. This remark was recently taken up by Dyson 
et al. \cite{dyson}, who proved that $|a_k/b_k-\pi|=O(Q^{-k})$ as 
$k\rightarrow\infty$ where 
$Q=\sqrt{1+(2\pi/\ln2)^2}$. Lehmer also showed that for positive integer $k$

\begin{equation}
S_k(z)=\frac{2^{k+z^{5/2}}z^{1/2}}{(4-z)^{k+{3/2}}}( A_k(z/4)
\sin^{-1}(\sqrt{z/4})+\sqrt{z(4-z)}B_k(z/4))
\label{2}
\end{equation}
where $A_k$ and $B_k$ are recursively defined  polynomials.  It was apparently 
not until 2005 that (\ref{2}) was evaluated explicitly, for $z=1$, by 
J. Borwein and P. Girgensohn \cite{bor}  who showed
\begin{equation}
S_k(1)=
\frac{1}{2}(-1)^{k+1}\sum_{j=1}^{k+1} (-1)^j j!S(j+1,j) 3^{-j} 
 \binom{2j}{j}\left(\sum_{i=1}^{j-1}\frac{3^i}{(2i+1) \binom{2i}{i}}
+ \frac{2}{3 \sqrt{3}}  \pi \right).
\label{3} 
\end{equation}
where the {\it Stirling numbers of the second kind} are defined by
\begin{equation}
S(k,j)=\frac{(-1)^j}{j!}\sum_{m=0}^j(-1)^m m^k\binom{j}{m}.
\label{4}
\end{equation}
The aim of the present  note is to extend (\ref{3}) to complex $z$ and thus to 
continue (\ref{1}) analytically beyond its circle of convergence $|z|=4$.

\section{Calculation}

We begin with the observation that  
$\left(m\binom{2m}{m}\right)^{-1} =B(m,m+1)$,
 where B denotes Euler's beta integral. Hence,
\begin{equation}
S_k(z)=\int_0^1\frac{dt}{t}\sum_{m=1}^{\infty}m^{k+1}(zt(1-t))^m.
\label{5}
\end{equation}

Next, equation (21) of Girgensohn and Borwein \cite{bor}, 
\begin{equation}
\sum_{m=1}^{\infty}m^pX^m=
\sum_{n=1}^p\sum_{m=1}^n(-1)^{m+n}\binom{n}{m} m^pX^n(1-X)^{-n-1},
\label{6}
\end{equation}
gives
\begin{equation}
S_k(z)=\sum_{n=1}^{k+1}\sum_{m=1}^n(-1)^{m+n}\binom{n}{m}
m^{k+1}\int_0^1\frac{dt}{t}\frac{(zt(1-t))^n}{(1-zt(1-t))^{n+1}}.
\label{7}
\end{equation}

In the appendix it is shown that
\begin{equation}
\int_0^1\frac{dt}{t}\frac{(zt(1-t))^n}{(1-zt(1-t))^{n+1}}=
\frac{\sqrt{\pi}\Gamma(n)}{\Gamma(n+1/2)}X^n\;_2F_1(-1/2,n;n+1/2;-X)
\label{8}
\end{equation}
\noindent
where $X=z/(4-z)$, so
\begin{equation}S_k(z)=
\sum_{n=1}^{k+1}n!B(n,1/2)S(k+1,n)X^n\;_2F_1(-1/2,n;n+1/2;-X).
\label{9}
\end{equation}
By induction, starting with the tabulated value for $n=1$ and using 
Gauss' contiguity relations we find (some details are given in the appendix)
\begin{equation}
_2F_1(-1/2,n;n+1/2;-X) = 
\nonumber
\end{equation}
\begin{eqnarray}
\left( \frac{1}{2} \right)_n \left(\frac{1}{n!}+
\frac{1}{\sqrt{\pi}\Gamma(n)}\sum_{k=0}^{n-1}\frac{(-1)^k\Gamma(k+1/2)}{(k+1)!}
\binom{n-1}{k}
\left(\frac{X+1}{X}\right)^{k+1}\times \right.  \nonumber \\
 \left.\left[\sqrt{X}\sin^{-1}\sqrt{\frac{X}{X+1}}-
\frac{1}{2}\sum_{l=1}^k\frac{(l-1)!}{(1/2)_l}
\left(\frac{X}{X+1} \right)^l \right] \right).
\label{10}
\end{eqnarray}
(We have used the ascending factorial notation $(a)_n=\Gamma(a+n)/\Gamma(a)$).
Therefore we have the principal result
\begin{equation}
S_k(z)=   
\sum_{n=1}^{k+1}n!\left(\frac{z}{4-z}\right)^nS(k+1,n)\times\nonumber\end{equation}
\begin{equation}
 \left(\frac{1}{n}+
\sum_{p=0}^{n-1}(-1)^p\frac{(1/2)_p}{(p+1)!}\binom{n-1}{p}
\left(\frac{4}{z}\right)^{p+1} \left(\sqrt{\frac{z}{4-z}}
\sin^{-1}\frac{\sqrt{z}}{2}-
\frac{1}{2}\sum_{l=1}^p\frac{\Gamma(l)}{(1/2)_l}\left(\frac{z}{4}\right)^l 
\right)  \right)
\label{11}
\end{equation}

Equation (\ref{11}) is rather condensed; in unpacking it, sums with upper 
limit less than the lower limit are to be interpreted as 0.  It is clear 
from (\ref{11})  that for {\it rational} $z$
\begin{equation}\sum_{m=1}^{\infty}\frac{n^kz^n}{\binom{2n}{n}}=
R_1(z,k)+R_2(z,k)\sqrt{\frac{z}{4-z}}\sin^{-1}\frac{\sqrt{z}}{2},
\label{12}
\end{equation}
where $R_j$ is a rational number. 

One sees from (\ref{11})  that $ S_k(z)$ is analytic on the two-sheeted 
Riemann surface formed by two planes cut and rejoined along the 
real half-line $x>4$. The numbers in (\ref{12}) have the explicit expressions
\begin{equation}
R_1(z,k)=\label{13} \end{equation}
\begin{equation}\sum_{n=1}^{k+1}n!S(k+1,n)\left(\frac{z}{4-z}\right)^n
\left(\frac{1}{n}-\frac{1}{2}\sum_{p=1}^{n-1}
\sum_{l=1}^p\frac{(-1)^p(1/2)_p}{(p+1)!(1/2)_l}\binom{n-1}{p}
\Gamma(l)\left(\frac{4}{z}\right)^{p-l+1}\right),
\nonumber
\end{equation}
\begin{equation}
R_2(z,k)=
\sum_{n=1}^{k+1}n!S(k+1,n)\sum_{p=0}^{n-1}\frac{(-1)^p}{(p+1)!}\binom{n-1}{p}
\left(\frac{4}{z}\right)^{p+1}.
\label{14}
\end{equation}

\section{Asymptotics}

It is convenient to work in terms of the exponential generating function
\begin{equation}
G(z,t):=\sum_{k=0}^{\infty}S_k(z)\frac{t^k}{k!}=S_0(ze^t)=
\frac{z}{4-ze^t}+\frac{4\sqrt{z}e^{t/2}}
{(4-ze^t)^{3/2}}\sin^{-1}\frac{\sqrt{z}e^{t/2}}{2}
\label{15}
\end{equation}
To find the generating functions  $\rho_j(z,t):=\sum R_j(z,k)t^k/k!$, it 
would be simplest to start with a series 
$D_k(z)=R_1(z,k)-R_2(z,k)\sqrt{\frac{z}{4-z}}\sin^{-1}\sqrt{z}/2$, work out 
its generating function $D(z,t)$ and by taking the sum and difference 
identify $\rho_1$ and $\rho_2$. However, this series has not been found and 
there is nothing to guarantee its existence in tractable form. Therefore, 
the $\rho_j$ were evaluated directly from (\ref{13}) and (\ref{14}). The 
details are omitted as the results
\begin{eqnarray}
\rho_1(z,t) & = & \frac{ze^t}{4-ze^t}+ \nonumber \\
& & \frac{8}{\pi}\sqrt{\frac{ze^t}
{(4-ze^t)^{3/2}}}\left(\sin^{-1}\frac{\sqrt{z}e^{t/2}}{2}
\cos^{-1}\frac{\sqrt{z}}{2}
-\cos^{-1}\frac{\sqrt{z}e^{t/2}}{2}\sin^{-1}\frac{\sqrt{z}}{2}\right),
\label{16}
\end{eqnarray}
\begin{equation}
\rho_2(z,t)=4\sqrt{\frac{(4-z)e^t}{(4-ze^t)^3}}
\label{17}
\end{equation}
are easily verified. In the case $z=2$, (\ref{15}) and (\ref{16}) are 
identical to Dyson's formulas \cite{dyson, dysonc} obtained empirically.

In view of the prominent role that the ratio $R_1(z,t)/R_2(z,t)$ plays in 
Dyson et al. \cite{dyson} for $z=2$ it is interesting to examine it 
for general $z$. From (\ref{17}) we have 
\begin{equation}
R_2(z,k)=
\frac{2k!\sqrt{4-z}}{\pi i}
\oint\frac{ds}{s^{k+1}}\frac{e^{s/2}}{(4-ze^s)^{3/2}}.
\label{18}
\end{equation}
The non-zero singularity closest to $s=0$ is $s_0=\ln(4/z)$ and it dominates 
the asymptotic behavior.
Ignoring the other singularities, distorting the contour to a small circle 
about $s_0$ and translating back to the origin by $t=s-s_0$, we have
\begin{equation}
R_2(z,k)\sim-\frac{k!\sqrt{4-z}}{zs_0^{k+1}}\oint\frac{dt}{2\pi i}
\frac{e^{t/2}}{(1-e^t))^{3/2}}.
\label{19}
\end{equation}
The exact value of the integral in (\ref{19}) is $-(2/\pi)\sqrt{e/(e-1)}$, and 
so 
\begin{equation}
R_2(z,k)\sim\frac{k!}{(\ln(4/z))^{k+1}}
\frac{2}{\pi}\sqrt{\frac{e(4-z)}{z(e-1)}}.
\label{20}
\end{equation}
In the same way we obtain
\begin{equation}
R_1(z,k)\sim
\frac{k!}{(\ln(4/z))^{k+1}}\left(\sqrt{2}+\frac{2}{\pi}\left(\sqrt\frac{e}{e-1}
-\sqrt{2}\right)\cos^{-1}\frac{\sqrt{z}}{2}-\frac{2^{3/2}}{\pi}
\sin^{-1}\frac{\sqrt{z}}{2}\right).
\label{21}
\end{equation}

\section{Discussion}

From (\ref{20}) and (\ref{21}) we find
\begin{equation}
\lim_{k\to\infty}\left(\frac{R_1(z,k)}{R_2(z,k)}-
\sqrt{\frac{z}{4-z}}\sin^{-1}\frac{\sqrt{z}}{2}\right)=
\sqrt{\frac{z}{4-z}}\left(\cos^{-1}\frac{\sqrt{z}}{2}-
\sin^{-1}\frac{\sqrt{z}}{2}\right).
\label{22}
\end{equation}
It thus appears that  Lehmer's choice, $z=2$, is the unique permissible case 
for which the limit (\ref{22}) vanishes. (Also the {\it Lehmer limit}, as 
defined by Dyson et al. \cite{dyson}, relates to $\pi/4$ here rather 
than $\pi$). Finally, for negative integer indices, since
\begin{equation}
2 S_{-k}(z)=\;_{k+1}F_k(1,\dots,1;\tfrac{3}{2},2,\dots,2;
\tfrac{1}{4}z),
\label{23}
\end{equation}
the fact that $S_{-k}(z)$ can be obtained from $S_0(z)$ by successive 
integrations with respect to $z$ and the explicit evaluations by 
Lehmer \cite{lehm}, Borwein and Girgensohn \cite{bor} and others 
\cite{6,7,8,9, 10} it should be possible to obtain  explicit values 
for sundry generalized hypergeometric functions.

\section{ Appendix: Derivation of Equations (\ref{8}) and (\ref{10})}

Let us consider, for any integrable function $F$,
\begin{equation}I=\int_0^1\frac{dt}{t}F(t(1-t))\nonumber\end{equation}
Let $u=t(1-t)$, so $u(0)=u(1)=0;\;\; u(1/2)=1/4$. Then there are two 
expressions for $t$:
\begin{equation}
t_+=\frac{1}{2}(1+\sqrt{1-4t})\quad\mbox { for }\frac{1}{2}
\le t\le1, \, \text{ with }
\frac{dt_+}{t_+}=\left(1-\frac{1}{\sqrt{1-4u}}\right)\frac{du}{u}
\nonumber\end{equation}
\noindent
and 
\begin{equation}
t_-=\frac{1}{2}(1-\sqrt{1-4t})\quad\mbox { for }0\le t\le\frac{1}{2}, 
\, \text{ with } 
\frac{dt_+}{t_+}=\left(1+\frac{1}{\sqrt{1-4u}}\right)\frac{du}{u}.
\nonumber
\end{equation}
Consequently, 
\begin{eqnarray}
I & = & \int_0^{1/2}\frac{dt_-}{t_-}F(u)+\int_{1/2}^1\frac{dt_+}{t_+}F(u)
=2\int_0^{1/4}\frac{du}{u\sqrt{1-4u}}F(u) \nonumber \\
& = & 
2\int_0^1\frac{dx}{x\sqrt{1-x}}F \left(\frac{1}{4}x \right)=
2 \int_0^1\frac{dt}{(1-t)\sqrt{t}}F\left(\frac{1-t}{4}\right)
\nonumber
\end{eqnarray}
and, with $t=x^2$,
\begin{equation}
I=4\int_0^1\frac{dx}{1-x^2}F\left(\frac{1-x^2}{4}\right).
\nonumber
\end{equation}

\vskip .1in
 
Therefore,
\begin{equation}
L=\int_0^1\frac{dt}{t}\frac{(zt(1-t))^{\alpha}}
{(1-zt(1-t))^{\beta}}=
2\left(\frac{z}{4}\right)^{\alpha-\beta}
\int_0^1dx\frac{(1-x^2)^{\alpha-1}}{(a^2+x^2)^{\beta}},
\nonumber
\end{equation}
where $a^2=1/X=(4-z)/z$.

From standard references
\begin{equation}
\int_0^1 dx \cos(xy)(1-x^2)^{\alpha-1}=
\sqrt{\frac{\pi y}{8}}\left(\frac{2}{y}\right)^{\alpha}
\Gamma(\alpha)J_{\alpha-1/2}(y),
\nonumber
\end{equation}
\begin{equation}
\int_0^{\infty}dx \cos(xy)(a^2+x^2)^{-\beta}
=\frac{\sqrt{\pi}}{\Gamma(\beta)}\left(\frac{y}{2a}\right)^{\beta-1/2}
K_{\beta-1/2}(ay)
\nonumber
\end{equation}
so, by the Parseval relation for the Fourier transform
\begin{equation}
\int_0^1dx\frac{(1-x^2)^{\alpha-1}}{(a^2+x^2)^{\beta}}=
\frac{2^{\alpha-\beta}}{a^{\beta-1/2}}\frac{\Gamma(\alpha)}{\Gamma(\beta)}
\int_0^{\infty}dy y^{\beta-\alpha}J_{\alpha-1/2}(y)K_{\beta-1/2}(ay).
\nonumber
\end{equation}
This is a tabulated Hankel Transform and yields
\begin{equation}L=
\sqrt{\pi}\left(\frac{4}{z}\right)^{\beta-\alpha}X^{\beta}
\frac{\Gamma(\alpha)}{\Gamma(\alpha+1)}\;
_2F_1 \left(\frac{1}{2},\beta;\alpha+\frac{1}{2};-X \right).
\nonumber
\end{equation}

Consequently
\begin{equation}
\int_0^1\frac{dx}{x}\frac{(zt(1-t))^n}{(1-zt(1-t))^{n+1}}
=\sqrt{\pi}\left(\frac{4}{z}\right)X^{n+1}\;
_2F_1 \left(\frac{1}{2},n+1;n+\frac{1}{2};-X \right)
\nonumber
\end{equation}
However, since $\;_2F_1(a,b;c;z)=(1-z)^{c-a-b}\;_2F_1(c-a,c-b;c;z)$,
\begin{equation}
\;_2F_1 \left(\tfrac{1}{2},n+1;n+\tfrac{1}{2};-X \right)=
(1+X)^{-1}\;_2F_1 \left(-\tfrac{1}{2};n;n+\tfrac{1}{2};-X \right)
\nonumber
\end{equation}

\vskip .1in

   Next, we note that \cite[p .\ 590]{prud}
   \begin{equation}
\;_2F_1(-1/2,1;3/2;z)=
\frac{1}{2} \left(1+(1-z)\frac{\tanh^{-1}\sqrt{z}}{\sqrt{z}} \right).
\nonumber
\end{equation}
   With $z\rightarrow-z$, noting that 
$-i \tanh^{-1} iw=\sin^{-1}\sqrt{\frac{w}{1+w}}$ one has
   \begin{equation}
\;_2F_1(-1/2,1;3/2;-z)=
\frac{1}{2}(1+(1+z)\frac{\sin^{-1}\sqrt{\frac{z}{1+z}}}{\sqrt{z}}).
\label{24}
\end{equation}
   We next apply Gauss' differentiation formula
   \begin{equation}
\frac{d}{dz}((1+z)^k\;_2F_1(-1/2,k;k+1/2;-z))=
\nonumber
\end{equation}
   \begin{equation}
\frac{2k(k+1)}{2k+1}(1+z)\;_2F_1(-1/2,k+1;k+3/2;-z).
\label{25}
\end{equation}
   Iteration of (\ref{25}) starting with (\ref{24}), after a great deal of 
tedious algebra, aided by Mathematica, results in (\ref{10}).

\section{Acknowledgements}

The author is grateful to Profs. N. E. Frankel for suggesting this 
problem and F. J. Dyson for enlightening correspondence.

\bigskip
\hrule
\bigskip

\noindent 2000 {\it Mathematics Subject Classification}:  Primary
11B65; Secondary 33B05.

\noindent {\it Keywords}: Binomial coefficient, infinite series, generating function.

\bigskip
\hrule
\bigskip

\noindent (Concerned with sequences 
\seqnum{A008277, A145557}

\bigskip
\hrule
\bigskip

\vspace*{+.1in}
\noindent
Received ; 
revised version received

\bigskip
\hrule
\bigskip

\noindent
Return to
\htmladdnormallink{Journal of Integer Sequences home page}{http://www.cs.uwaterloo.ca/journals/JIS/}.
\vskip .1in

\end{document}